\documentclass{amsart}
\usepackage{amsmath}
\usepackage{amscd}
\usepackage{amssymb}
\usepackage{amsthm}
\RequirePackage{filecontents}
\usepackage{fullpage}
\usepackage[numbers]{natbib}  
\usepackage[draft]{hyperref}

\theoremstyle{plain}
\newtheorem{prop}{Proposition}

\theoremstyle{definition}

\newtheorem{rem}[prop]{Remark}
\newtheorem{ex}[prop]{Example}

\author{Brandon Williams }

\address{Department of Mathematics \\ University of California \\ Berkeley, CA 94720}
\email{btw@math.berkeley.edu}

\begin{document}

\nocite{*}

\begin{abstract} This note points out that for any odd prime $p \ge 5$, Zagier's weight $3/2$ mock Eisenstein series can be completed to a $p$-adic modular form in a way that bears some resemblance to its completion to a harmonic Maass form. \end{abstract}

\title{A $p$-adic completion of Zagier's Eisenstein series}

\maketitle

The Eisenstein series of weight two, $E_2(\tau) = 1 - 24 \sum_{n=1}^{\infty} \sigma_1(n) q^n,$ is not a modular form. However, the coefficients of $E_2$ are congruent to those of $E_{p+1}$ mod $p$ for every odd prime $p$, making $E_2$ a ``mod $p$ modular form'' as in \cite{S}. More precisely, $E_2$ is a $p$-adic modular form in Serre's sense; that is, it is a $p$-adic limit of modular forms with integer coefficients (and possibly varying weight). \\

It is natural to look for similar results among more complex Eisenstein series. One family of candidates are the \textbf{Cohen Eisenstein series} introduced in \cite{C}, which have half-integral weight and satisfy Kohnen's plus space condition. For $k \in 3/2 + 2 \mathbb{Z}$, $k \ge 7/2$ the Cohen Eisenstein series is \begin{equation} E_k(\tau) = \sum_{n=0}^{\infty} c_{n,k} q^n = 1 + \frac{1}{\zeta(2 - 2k)} \sum_{n \equiv 0,3 \, (4)} \Big[ L(3/2 - k, \chi_{-D_0}) \sum_{d | f} \mu(d) \chi_{-D_0}(d) d^{k-3/2} \sigma_{2k-2}(f/d) \Big] q^n,\end{equation} where we write $n = D_0 f^2$ with $-D_0$ being the discriminant of $\mathbb{Q}(\sqrt{-n})$, and $\chi_{-D_0}(m) = \left( \frac{-D_0}{m} \right)$ is the associated Dirichlet character with $L$-function $L(s,\chi_{-D_0}) = \sum_{n=1}^{\infty} \chi_{-D_0}(n) n^{-s}$. The smallest examples are \begin{align*} E_{7/2}(\tau) &= 1 + 56q^3 + 126q^4 + 576q^7 + 756q^8 + 1512q^{11} + 2072q^{12} + ... \\ E_{11/2}(\tau) &= 1 - 88q^3 - 330q^4 - 4224q^7 - 7524q^8 - 30600q^{11} - 46552q^{12} - ... \\ E_{15/2}(\tau) &= 1 + 56q^3 + 366q^4 + 14016q^7 + 33156q^8 + 260712q^{11} + 462392q^{12} + ... \end{align*} The easiest way to explain the half-integral weight transformation law is that $E_k(\tau) \theta(\tau)$ is a modular form of weight $k+1/2 \in 2 \mathbb{Z}$ and level $\Gamma_0(4)$, where $\theta(\tau) = 1 + 2q + 2q^4 + 2q^9 + ...$ is the usual theta series. \\

Setting $k = 3/2$ in the coefficient formula (1) produces Zagier's mock Eisenstein series $$E_{3/2}(\tau) = 1 - 4q^3 - 6q^4 - 12q^7 - 12q^8 - 12q^{11} - 16q^{12} - ...$$ in which the coefficient of $q^n$ is $-12$ times the Hurwitz class number $H(n)$. (Elsewhere this is often rescaled to have constant term $-1/12$.) For $k \ge 7/2$, congruences such as $E_k \equiv E_{k+p-1} \, (\text{mod} \, p)$ continue to hold for the Cohen Eisenstein series (as observed by Koblitz, \cite{K}). However, the congruence $E_{3/2} \equiv E_{3/2 + p - 1}$ mod $p$ does not hold. For example, with $p = 5$ the coefficients of $q^3$ in $E_{3/2}$ and $E_{11/2}$ are not congruent. Instead we can prove the following remark:

\begin{prop} Let $p$ be an odd prime. There are $p$-adic integers $a_m$, which are zero unless $m$ has the form $m \equiv -n^2 \, \text{mod} \, p$ for some $n \in \mathbb{Z}$, such that $$E_{3/2}(\tau) + \sum_{m=0}^{\infty} a_m q^m$$ is a $p$-adic modular form in Serre's sense.
\end{prop}

This completion almost resembles the completion of $E_{3/2}$ to a harmonic Maass form \cite{Z}, which also changes only the Fourier coefficients in exponents of the form $-n^2$: $$E_{3/2}^*(\tau) = E_{3/2}(\tau) - \frac{12}{\sqrt{y}} \sum_{n=-\infty}^{\infty} \beta(4\pi n^2 y) q^{-n^2},$$ where $\beta(x) = \frac{1}{16\pi} \int_1^{\infty} u^{-3/2} e^{-xu} \, \mathrm{d}u$. This observation may not be new but it does not seem to be readily available in the literature. It would be very interesting to find some completion of $E_{3/2}$ that accounts for both of these results at once, but I do not know any way to do this. In any case, this seems to be at least an interesting coincidence. \\

\begin{proof} To be explicit we write $$E_{3/2}(\tau) = 1 - 12 \sum_{n \equiv 0,3 \, (4)} \Big[ L(0,\chi_{-D_0}) \sum_{d | f} \mu(d) \chi_{-D_0}(d) \sigma_1(f/d) \Big].$$ For $l \in \mathbb{N}$ consider the Cohen Eisenstein series of weight $k = 3/2 + p^{l-1} (p-1)$. The congruences $$\sigma_{2k-2}(f/d) \equiv \sigma_1(f/d), \; \; d^{k-3/2} \equiv \begin{cases} 1: & p \nmid d; \\ 0: & p | d \end{cases}$$ mod $p^l$ are valid for elementary reasons, and the inverse zeta value is $$\zeta(2 - 2k)^{-1} \equiv \zeta(-1)^{-1} \equiv -12 \, \bmod \, p;$$ in particular, reducing $\frac{1}{\zeta(2-2k)}$ mod $p^l$ is well-defined even when $p=3$. To handle the $L$-function values we use Kubota and Leopoldt's $p$-adic $L$-functions $L_p(s,\chi)$ attached to Dirichlet characters $\chi$, which have the special values $$L_p(1-n,\chi \omega^n) = (1 - \chi(p)  p^{n-1}) L(1-n,\chi), \; \; n \in \mathbb{N},$$ where $\omega$ is the Teichm\"uller character mod $p$. In particular, for $n > 1$ the congruence $$L_p(1-n,\chi \omega^n) \equiv L(1-n,\chi)$$ holds mod $p^{n-1}$, while for $n=1$ we get exact values $$L_p(0,\chi \omega) = \begin{cases} 2L(0,\chi): & \chi(p) = -1; \\ L(0,\chi): & \chi(p) = 0; \\ 0: & \chi(p) = 1. \end{cases}$$ (See \cite{W}, chapter $5$ for details.) Kummer's congruences imply that $L_p(1-m,\chi) \equiv L_p(1-n,\chi) \, \text{mod} \, p^l$ whenever $m \equiv n \, \text{mod} \, p^{l-1}$, so $$L(0,\chi_{-D_0}) \equiv \frac{1}{2}L_p(0,\chi \omega) \equiv \frac{1}{2} L_p((1-p)p^{l-1}, \chi \omega) = \frac{1}{2} L_p((1-p)p^{l-1}, \chi \omega^{1 + p^{l - 1}(p-1)}) \equiv \frac{1}{2}L((1-p)p^{l-1}, \chi) \, (\text{mod} \, p^l) $$ for all discriminants $-D_0$ with $\chi_{-D_0}(p) = -1$. In other words, the coefficients $c_{m,k}$ of $E_k$ for this particular weight satisfy $$c_{m,k} \equiv \frac{1 - \chi_{-m}(p)}{\zeta(2-2k)} H(m) \; \text{mod}\, p^l.$$ Using the $p$-adic limit $$\lim_{l \rightarrow \infty} \zeta(-1 - p^l (p-1)) = (1-p) \zeta(-1) = \frac{p-1}{12} \in \mathbb{Q}_p,$$ we see that $E_{3/2}(\tau) = -12\sum_m H(m) q^m$ differs from the $p$-adic modular form $$\lim_{l \rightarrow \infty} -6 \zeta(-1 - p^l (p-1))E_{3/2 + p^l (p-1)}(\tau) \equiv \frac{1-p}{2} \lim_{l \rightarrow \infty} E_{3/2 + p^l (p-1)} $$ only by coefficients of $q^m$ for which $\chi_{-m}(p) \ne 1$, or equivalently for $m$ not of the form $-n^2$ mod $p$.
\end{proof}

\begin{rem}
This argument makes clear how one can choose the numbers $a_m$: namely, we can take $a_0 = -\frac{1+p}{2}$ and for $m > 0$, $$a_m = (p-1) H(m)$$ if $m \equiv 0$ mod $p$ and $$a_m = 2 (p-1) H(m)$$ if $\left( \frac{m}{p} \right) = 1.$ (Note that $a_m$ is not always integral for primes $p \equiv 2 \, (\text{mod} \, 3)$.)
\end{rem}

\begin{ex} For the prime $p = 7$, reducing mod $p$ gives the series \begin{align*} E_{3/2}(\tau) &\equiv 1 + 3q^3 + 1q^4 + 2q^7 + 2q^8 +2q^{11} + 5q^{12} + 4q^{15} + 3q^{16} + ... \, (\text{mod} \, 7), \\ \frac{1-7}{2} E_{15/2}(\tau) &\equiv 4 + 0q^3 + 1q^4 + 1q^7 + 2q^8 + 2q^{11} + 0q^{12} + 4q^{15} + 3q^{16} + ... \, (\text{mod} \, 7), \end{align*} which differ in the coefficients of $q^m$ with $$m = 0 \equiv -0^2, \; m = 3 \equiv -2^2, \; m = 7 \equiv -0^2, \; m = 12 \equiv -3^2,... \, (\text{mod} \, 7).$$ Modulo $7^2$, we should compare \begin{align*} E_{3/2}(\tau) &\equiv 1 + 45q^3 + 43q^4 + 37q^7 + 37q^8 + 37q^{11} + 33q^{12} + 25q^{15} + 31q^{16} + ... \, (\text{mod} \, 49), \\ \frac{1-7}{2} E_{87/2}(\tau) &\equiv 46 + 0q^3 + 43q^4 + 43q^7 + 37q^8 + 37q^{11} + 0q^{12} + 25q^{15} + 31q^{16} + ... \, (\text{mod} \, 49). \end{align*}
\end{ex}

\noindent \textbf{Acknowledgements:} I would like to thank Ken Ribet for reading this note and providing helpful feedback.

\bibliographystyle{plainnat}
\bibliography{\jobname}
\end{document}